\documentclass{article}
 
\usepackage{amssymb}
\usepackage{color}
\usepackage{graphicx}
\usepackage{enumerate}

\begin{document}
\newtheorem{Theorem}{Theorem}[section]
\newtheorem{Proposition}[Theorem]{Proposition}
\newtheorem{Lemma}[Theorem]{Lemma}
\newtheorem{Example}[Theorem]{Example}
\newtheorem{Corollary}[Theorem]{Corollary}
\newtheorem{Fact}[Theorem]{Fact}
\newtheorem{Conjecture}[Theorem]{Conjecture}
\newenvironment{Definition} {\refstepcounter{Theorem} \medskip\noindent
 {\bf Definition \arabic{section}.\arabic{Theorem}}\ }{\hfill}
\newenvironment{Remarks} {\refstepcounter{Theorem}
\medskip\noindent {\bf Remarks
\arabic{section}.\arabic{Theorem}}\ }{\hfill}

\newenvironment{Question} {\medskip\refstepcounter{Theorem}
     \noindent {\bf Question
\arabic{section}.\arabic{Theorem}}\ } {\hfill}

 \newcommand{\qed}{\hfill \ensuremath{\Box}}
\newenvironment{Proof}{{\noindent \bf Proof\ }}{\hfill\qed}

\newenvironment{claim} {{\smallskip\noindent \bf Claim\ }}{\hfill}

\def \blue {\color{blue}}
\def\red{\color{red}}
\def \L {{\cal L}}

\def \B {{\cal B}}
\def \mod {{\rm mod \ }}
\def \iso {\cong}
\def \Lor  {\L_{\rm or}}
\def \Lr {\L_{\rm r}}
\def \Lg {\L_{\rm g}}
\def \I {{\cal I}}
\def \M {{\cal M}}\def\N {{\cal N}}
\def \E {{\cal E}}
\def \Proj {{\mathbb P}}
\def \H {{\cal H}}
\def \x {\times}
\def \Stab {{\rm Stab}}
\def \Z {{\mathbb Z}}
\def \V {{\mathbb V}}
\def \C {{\mathbb C}} \def \Cexp {\C_{\rm exp}}
\def \R {{\mathbb R}}
\def \Q {{\mathbb Q}}\def \K {{\mathbb K}}
\def \F {{\cal F}}\def \A {{\mathbb A}}
 \def \X {{\mathbb X}}
\def \G {{\mathbb G}}
\def\HH {{\mathbb H}}
\def \Nn {{\mathbb N}}\def \Nn {{\mathbb N}}
\def\D {{\mathbb D}}
\def \hat {\widehat}
\def \bar{\overline}
\def \Spec {{\rm Spec}}
\def \bul {$\bullet$\ }
\def\proves {\vdash}
\def \Co {{\cal C}}
\def \ACFp {{\rm ACF}_p}
\def \ACF0 {{\rm ACF}_0}
\def \ee {\prec}
\def \Diag {{\rm Diag}}
\def \Diage {{\rm Diag}_{\rm el}}
\def \DLO {{\rm DLO}}
\def \d {{\rm depth}}
 \def \dist {{\rm dist}}
\def \P {{\cal P}}
\def \ds {\displaystyle}
\def \Fp {{\mathbb F}_p}
\def \acl {{\rm acl}}
\def \dcl {{\rm dcl}}

\def \dom {{\rm dom\ }}
\def \tp {{\rm tp}}
\def \stp {{\rm stp}}
\def \Th  {{\rm Th}}
\def\< {\Lngle}
\def \> {\rangle}
\def \n {\noindent}
\def \minusdot{\hbox{\ {$-$} \kern -.86em\raise .2em \hbox{$\cdot \
$}}}
\def\exp {{\rm exp}}\def\ex {{\rm ex}}
\def \td {{\rm td}\ }
\def \ld {{\rm ld}}
\def \span {{\rm span}}
\def \tilde {\widetilde}
\def \d {\partial}
\def \del {\partial}
\def \cl {{\rm cl}}
\def \acl {{\rm acl}}
\def \cN {{\cal N}}
\def \Qalg {{\Q^{\rm alg}}}
\def \th {^{\rm th}}
\def \deg { {\rm deg} }
\def\hat {\widehat}
\def\li {\L_{\infty,\omega}}
\def\lo {\L_{\omega_1,\omega}}
\def\lk {\L_{\kappa,\omega}}
\def \ee {\prec}\def \bSigma {{\mathbf\Sigma}}
\def \mod {\ {\rm mod\ }}
\def \Tor {{\rm Tor}}
\def \td {{\rm td}}
\def \dim {{\rm dim\ }}
 \def \| {\kern -.3em \restriction \kern -.3em}
 \def \lc{\lceil}
 \def \rc{\rceil}
 \def \SR {{\rm SR}}
 \def\L {{\cal L}}
 \def \a{{\bf a}}
 \def \b{{\bf b}}
 \def\x{{\bf x}}
 \def \< { \langle}
 \def \> { \rangle}
 \def \K {{\cal K}}
  \def \MM {{\mathbb M}} 
  \def \vm {{\mathfrak m}}
  \def\img{{\rm img\ }}
  \def\calR {{\cal R}}

 \title{Rigid Real Closed Fields}
 \author{David Marker\\ University of Illinois Chicago\and Charles Steinhorn\\ Vassar College 
\thanks{Partially supported by Simons Foundation Mathematics and Physical Sciences Collaboration
Grant for Mathematicians~524012.}  
}
\date{}

\maketitle
 
 \begin{abstract} We construct a non-Archimedean real closed field of transcendence degree two with no nontrivial automorphisms. 
 This is the first construction of a countable rigid non-Archimedean real closed field,
 \end{abstract}
 
We say that a structure is {\em rigid} if its   automorphism group is trivial. In real closed fields, the positive elements are the nonzero squares. Thus
every automorphism of a real closed field preserves the ordering as well as the field structure.  If $K$ is an Archimedean real closed field, the field of rational numbers is dense in $K$ and fixed pointwise by all automorphisms; hence $K$ is rigid.  

Are there non-Archimedean rigid real closed fields?  In \cite{ss},  Shelah proved that it is consistent with ZFC that there are.  More specifically, 
he showed, assuming Jensen's combinatorial principle $\diamondsuit_{\kappa^+}$,  that there are rigid non-Archimedean real closed fields of cardinality
$\kappa^+$.   In later work,   Mekler and Shelah \cite{ms} revisited Shelah's results and showed that the existence of arbitrarily large rigid non-Archimedean real closed fields
could be proved in ZFC without extra set-theoretic assumptions.  We are grateful to Biran Falk Dotan for calling our attention to the Mekler--Shelah paper. 

In 2018 Ali Enayat asked on MathOverflow if there are countable rigid non-Archimedean real closed fields.  Our main result, Theorem \ref {theorem}, is that there are rigid non-Archimedean real closed fields of transcendence degree two.    This gives the first positive answer to Enayat's question.

\section{Preliminaries}

 Let $k$ be the field of real algebraic numbers. We work in a sufficiently saturated real closed field $\calR\supset\R $.  All the fields we construct are subfields of $\cal R$. For a real closed field $K$, let
$K\< a_1,\dots a_n\> $ denote the real closure of the ordered field $K(a_1,\dots,a_n)$. The rigid field we construct in Theorem \ref {theorem} has the form $k\< a, b\> $, where $a$ is an infinite element and $b$ is transcendental over $k\< a\> $.

If $K\subset \calR$ is real closed and $a_1,\dots,a_n\in \calR$, we let $\tp(a_1,\dots,a_n/K)$, the {\em type} of $a_1,\dots,a_n$ over $K$,
be the set of all first order formulas $\phi(v_1,\dots, v_n)$ in the language of ordered rings with free variable from $v_1,\dots, v_n$ and parameters from $K$ such that
${\calR}\models \phi(a_1,\dots,a_n)$.  By Tarski's quantifier elimination, $\tp(a_1,\dots,a_n/K)$ is completely determined by knowing   
$$\{p\in K[X_1,\dots,X_n] : p(a_1,\dots, a_n)\ge 0\}.$$ When $n=1$ this simplifies further, and $\tp(a/K)$ is completely determined by
$$\{v\le m: m\in K, a\le m\}\cup\{v\ge m: m\in K, a\ge m\},$$  
i.e., the position of $a$ with respect to the ordering of $K$.
 
 We use several basic facts about definable sets and functions in real closed fields that hold in every o-minimal expansion of a real closed field (see
 \cite{ps} and \cite{vdd} for details):   
\begin{itemize} 
\item (o-minimality) If $K\subset\calR$ is a real closed, every $K$-definable subset of $\calR$ is a finite union of points and intervals with endpoints in $K\cup\{\pm \infty\}$.

\item(monotonicity) If $I\subset  \calR$ is an interval, and $f:I\rightarrow\calR$ is definable in $\calR$, we can partition $I=J_1\cup\dots\cup J_m\cup X$,
such that $X$ is finite and each $J_i$, for $i\leq m$, is an (open) interval with endpoints in $\calR \cup\{\pm \infty\}$ on which $f\restriction J_i$ is continuous and strictly monotonic. 

\item (exchange) If $K\subset\calR$ is real closed, $a,b\in\cal R$ and $b\in K\< a\> $, then $b\in K$ or $a\in K\< b \> $.

\item (algebraic closure = definable closure) If $K\subset\calR$ is real closed, $A\subset \calR$, and $b\in K\< A\> $, then there is a $K$-definable function 
$F:\calR^n\rightarrow \calR$ and $a_1,\dots,a_n\in A$ such that $b=F(a_1,\dots,a_n)$.

\item (cell decomposition) For $K\subset\calR$ real closed, we say that $C\subseteq \calR^2$ is a {\em cell} defined over $K$ if either:

\begin{itemize} 
\item $C$ is a point in $K^2$;
\item $C=\{a\}\times I$ where $a\in K$ and $I$ is a $K$-definable interval;

\item there is an interval $I$ and a $K$-definable continuous $f:I\rightarrow\calR$ such that $C$ is the graph of $f$; 

\item there is an interval $I$ and $K$-definable functions $g_0,g_1:I\rightarrow \calR$ such that $g_0(x)<g_1(x)$ for all $x\in I$ and 
$$C=\{(x,y): x\in I, g_0(x)<y<g_1(x)\}.\footnote{We allow the possibilities that $g_0$ is identically $-\infty$  or $g_1$ is identically $+\infty$, but cells of this form will not be relevant in our construction.}$$
\end{itemize}
The Cell Decomposition Theorem asserts that every $K$-definable subset of $\calR^2$ is a finite union of cells defined over $K$.
 \end{itemize}

\section{Construction of a rigid non-Archimedean real closed field}

We begin by noting that we cannot improve our main result to transcendence degree one.  Recall that $k$ is the field of real algebraic numbers.

\begin{Proposition}\label{td1} If $K$ is a  non-Archimedean real closed field of transcendence degree one, then $K$ has a nontrivial automorphism.  Indeed,
$|{\rm Aut}(K)|=\aleph_0$.
\end{Proposition}

\begin{Proof} Let $a\in K$ with $a$ infinite.  By exchange, $K=k\< a\> $. Let $b\in k\< a\> $ be any other infinite element. For example, $b=a^m$ for some $m>1$. Since $a$ and $b$
realize the same cut over $k$, we have $\tp(a/k)=\tp(b/k)$.  Thus there is an ordered field isomorphism $\sigma$ between $k(a)$ and $k(b)$ where $\sigma|k$ is the identity and $\sigma(a)=b$.
By the uniqueness of 
real closures, $\sigma$ extends to an isomorphism from $K$ onto $k\< b\> $.  By exchange, $K =k\< b\> $ and there is an automorphism
of $K$ sending $a$ to $b$.

Moreover, if $\sigma\in {\rm Aut}(K)$, then, because definable closure and algebraic closure agree, $\sigma$ is completely determined by knowing $\sigma(a)$.
Thus $|{\rm Aut}(K)|=\aleph_0$.   
\end{Proof}

\begin{Theorem}\label{theorem} There is a rigid non-Archimedean real closed field of transcendence degree two.
\end{Theorem}

The field we construct has the form $K=k\< a, b\> $ where $a$ is infinite and $b$ is transcendental over $k\< a\> $.  We first argue that to make $K$ rigid it is both necessary and sufficient to show that $(a,b)$ is the unique realization of $\tp(a,b/k)$
in $K$.

To see necessity, let $u,v \in k\< a, b\> $
and $\tp(a,b/k)=tp(u,v/k)$.  Then, arguing as in the proof of Proposition~\ref{td1}, the ordered fields $k(a,b)$ and $k(u,v)$ are isomorphic
and this isomorphism extends to their real closures $K$ and $k\< u, v\> $, respectively. By exchange $k\< u, v\> =K$ and there is a 
nontrivial automorphism of $K$. For sufficiency, note that if $k\< a, b\> $ is not rigid, it must have an automorphism $\sigma$ such that $\sigma (a,b)=(u,v)\not=(a,b)$.

Thus, to build a rigid $K$ we need to ensure that $(a,b)$ is the only realization of $\tp(a,b/k)$ in $K$. Observe that because algebraic closure agrees with definable closure, we must construct $\tp(a,b/k)$ such that if $F:\calR^2\rightarrow \calR^2$ is definable over
$k$ and $F(a,b)\ne (a,b)$, then $\tp(F(a,b)/k)\ne\tp (a,b/k)$. 
As we want $a>k$ and $b$ transcendental over $k\< a\> $, there are restrictions on the $k$-definable cells $C$ where $(a,b)\in C$. 

\begin{Definition} We say that $C\subset \calR^2$ is an {\em end-cell} if there is $\alpha\in k$ and $k$-definable and continuous $h_0,h_1:(\alpha,+\infty)\rightarrow\calR$ such that $h_0(x)<h_1(x)$ for all $x>\alpha$ and 
$$C=\{(x,y): x>\alpha, h_0(x)<y<h_1(x)\}.$$
\end{Definition}

Observe that if $C$ is an end-cell and $X\subseteq C$ is $k$-definable, then, from a cell decomposition of $X$, we see that  there is an end-cell 
$C^\prime\subseteq C$ such that either $C^\prime\subseteq X$ or
$C^\prime\cap X=\emptyset$.

We claim that $\tp(a,b/k)$ is determined by the end-cells that contain $(a,b)$.
Indeed, suppose $C\subseteq\calR^2$ is a cell defined over $k$ and $(a,b)\in C$. 
As we want $a>k$, the $x$-coordinates of $C$ cannot be bounded.  
Also, $C$ must have dimension 2, as otherwise $b$ would be algebraic over $k(a)$.
Moreover, were $|b|>k\< a\> $, as in the proof of Proposition~\ref{td1}, we could build non-trivial automorphisms of $K$ fixing $k\< a\> $. 
Thus we can find $\alpha\in k$ and continuous, $k$-definable $h_0,h_1:(\alpha,+\infty)\rightarrow\calR$ 
with $h_0(x)<h_1(x)$ for all $x>\alpha$ such that 
$$(a,b)\in \{(x,y):x>\alpha, h_0(x)<y<h_1(x)\}\subseteq C.$$
 
The principal ingredient needed for the proof of Theorem~\ref{theorem} is 

\begin{Lemma}\label{main} Suppose $C$ is an end-cell and $F:C\rightarrow \calR^2$ is $k$-definable. Then there is an end-cell $C^\prime\subseteq C$
such that  either $F\restriction C^\prime$ is the identity function or $\img (F\restriction C^\prime)\cap C^\prime=\emptyset.$\end{Lemma}

\n{\bf Proof of Theorem~\ref{main} from Lemma~\ref{main}}
 
Let $F_0,F_1,\dots,F_n,\dots$ enumerate all $k$-definable functions $F:\calR^2\rightarrow\calR^2$ and let $\phi_0(v_1,v_2),\phi_1(v_1,v_2),\dots$
enumerate all formulas in the language of ordered rings with parameters from $k$ allowing only $v_1$ and $v_2$ as free variables.

Let $C_0=(0,+\infty)\times \calR$. 
We build a sequence of $k$-definable end-cells $$C_0\supset C_1\supset\dots\supset C_{n}\supset\dots.$$ 
Applying Lemma~\ref{main} we find an end-cell $\hat C_{n+1}\subset C_n$ such that either $F\restriction \hat C_{n+1}$ is the identity function or  
$\img (F_n\restriction \hat C_{n+1})\cap \hat C_{n+1}=\emptyset$. 
 We then can find an end-cell $C_{n+1}\subseteq \hat C_{n+1}$ such that $C_{n+1}\subseteq (n,+\infty)\times \calR$ and either $\phi_n(x,y)$ for all $(x,y)\in C_{n+1}$
 or $\neg\phi(x,y)$ for all $(x,y)\in C_{n+1}$.
 
We see that $\bigcap C_n$ determines a unique 2-type $p$ over $k$. If $(a,b)$ realizes $p$, then $a>k$ and if $F:\calR^2\rightarrow\calR^2$ is $k$-definable and $F(a,b)\ne (a,b)$, then $F(a,b)$ does not realize $p$. Thus $k\< a,b\> $ is non-Archimedean and has no nontrivial automorphisms.
 \hfill\qed
 
 \medskip
We have a great deal of freedom in our construction. It easily could be modified to build $2^{\aleph_0}$ non-isomorphic rigid non-Archimedean real closed fields of transcendence degree two.

 For the proof of Lemma~\ref{main}, we need the following well known multivariable analog of the monotonicity theorem.

\begin{Fact}\label{smooth cd}  If $C\subseteq\calR^2$ is an end-cell and $F:C\rightarrow\calR^2$ is $k$-definable, we can find a $k$-definable end-cell $C^\prime\subset C$ such that  either $F~\restriction C^\prime$ is a  continuous injection or
$\img (F\restriction C^\prime)$ has dimension at most 1.
\end{Fact}

 Fact \ref{smooth cd} is an easy consequence of trivialization in o-minimal theories (\cite{vdd} 7.1.2).  A more elementary proof could be given using cell decomposition, uniform finiteness and results on the dimension of fibers of definable maps (\cite{vdd} 4.1.6).

\medskip
\n{\bf Proof of Lemma \ref{main}}

Let $C$ be the end-cell $C:=\{(x,y): x>\alpha, h_0(x)<y<h_1(x)\}$. There are several cases to consider. 
 
 \medskip \n {\bf Case 1}\enspace There is an end-cell $C_0\subseteq C$ such that $\img (F\restriction C_0)$ has dimension at most one.
 
 In this case we can find $C^\prime\subseteq C_0$ such that $\img (F\restriction C^\prime)\cap C^\prime=\emptyset$, yielding the conclusion of the lemma.
 
 \medskip If we are not in Case~1, then, by Fact~\ref{smooth cd}, without loss of generality we may assume that $F\restriction C$ is a continuous injection.

\medskip\n {\bf Case 2}\enspace There is an end cell $C^\prime\subseteq C$ such that $F\restriction C^\prime$ is the identity function.

In this case the conclusion of the lemma certainly holds. 

\medskip If we are not in Case~2, then, by cell decomposition we may assume 
without loss of generality that $F(x,y)\ne (x,y)$ for all $(x,y)\in C$. We introduce some notation for the remaining cases. 
Let $f:(\alpha,+\infty)\rightarrow \calR$ be $k$-definable such that $h_0(x)<f(x)<h_1(x)$ for all $x>\alpha$. 
We define $\mu_f, \nu_f:(\alpha,+\infty)\rightarrow\calR$ by $F(x,f(x))=(\mu_f(x), \nu_f(x)).$
 
\medskip
\n{\bf Case 3}\enspace $\displaystyle \lim_{x\to \infty} \mu_f(x)\ne +\infty$ for some such function $f$.

 By monotonicity,  there are $\alpha^\prime>\alpha$ and $\beta$ such that $\mu_f(x)<\beta$ for all $x>\alpha^\prime$ and $\alpha^\prime>\beta$. 
 Applying cell decomposition we can find $k$-definable  $g_0$ and $g_1$ such that $h_0(x)\le g_0(x)<f(x)<g_1(x)\le h_1(x)$ and, with $\pi$ denoting 
 projection onto the first coordinate,  $\pi(F(x,y))<\beta$ for all 
$$(x,y)\in C^\prime:=\{(x,y): x>\alpha^\prime, g_0(x)<y<g_1(x)\}.$$ 
Then $\img (F\restriction C^\prime)\cap C^\prime=\emptyset $ and thus the conclusion of the lemma holds. 

\medskip Assuming we are not in Case~3, for every such function $f$ as above, $\displaystyle \lim_{x\to \infty} \mu_f(x)= +\infty$. By monotonicity,
increasing $\alpha$ if necessary, we may assume that $\mu_f$ is increasing.\footnote{At several points in the proof we shrink the cell by increasing $\alpha$.  We do this only finitely many times so always have an end-cell remaining.}Thus, for sufficiently large $x$, we can define 
$$f^*(x)=\nu_f(\mu^{-1}_f(x)).$$
The graph of $f^*$ is contained in the image, for sufficiently large $x$, of the graph of $f$ under $F$.  Increasing $\alpha$ if necessary,
we may assume that both $f$ and $f^*$ are defined on $(\alpha,+\infty)$ and that $F(x,f(x))$ is the graph of $f^*$ for sufficiently large $x$.
 
\medskip \n{\bf Case 4}\enspace There is a definable  function $f$ such that $f(x)\ne f^*(x)$ for sufficiently large $x$.

 Increasing $\alpha$ if necessary, we can assume that either $f^*(x)>f(x)$ for all sufficiently large $x$ or $f^*(x)<f(x)$ for all sufficiently large $x$.
The cases are similar and we consider only the alternative $f^*(x)>f(x)$ for sufficiently large $x$.

We claim we can find disjoint tubular neighborhoods of the graphs of $f$ and $f^*$ such that the neighborhood of the graph of $f$ is mapped by $F$ into the neighborhood of the graph of $f^*$.
 More precisely, we assert that we can find 
 $k$-definable 
 functions $g_0, g_1, \phi_0,\phi_1$ with domain $(\gamma,+\infty)$, where $\gamma\ge \alpha$, such that 
 $$h_0(x)<g_0(x)<f(x)<g_1(x)<h_1(x)\hbox{ and } g_1(x)< \phi_0(x)<f^*(x)<\phi_1(x)$$ 
 for all $x\in (\gamma,+\infty)$ and if we put  
 $$C^\prime:=\{(x,y): x>\gamma, g_0(x)<y<g_1(x)\}$$ 
 then $$\img (F\restriction C^\prime)\subseteq  \{(x,y):x>\gamma, \phi_0(x)<y<\phi_1(x)\}.$$
 In this case $\img (F\restriction C^\prime)\cap C^\prime=\emptyset$  and we have satisfied the conclusion of the lemma.
 
 To find the desired tubular neighborhood, first choose continuous $k$-definable $\phi_0$ and $\phi_1$ and $\beta$ such that 
 $$f(x)<\phi_0(x)<f^*(x)< \phi_1(x)$$ for all $x>\beta$.
 By continuity of $F$ and cell decomposition, we can find $k$-definable $g_0$ and $\hat g$ such that $g_0(x)<f(x)<\hat g(x)$ for all $x>\beta$ and 
 $\{(x,y): x>\beta \hbox{ and } g_0(x)<y<\hat g(x) \}$ is contained in $$F^{-1} (\{(x,y): x>\beta \hbox{ and } \phi_0(x)<y< \phi_1(x)\}.$$
Now define $$g_1(x)= f(x)+ {\min(\hat g(x), \phi_0(x))-f(x)\over 2}$$ and if necessary choose $\gamma\geq \beta$ so that $g_0$ and $g_1$ are continuous for $x>\gamma$.

 \medskip  \n{\bf Case 5}\enspace None of Cases~1-4 holds. 
 
 We complete the proof by showing that this case leads to a contradiction.
  
 For each $r\in [0,1]$ let $$f_r(x)= h_0(x)+ r(h_1(x)-h_0(x)).$$ We may assume that for each $r$ there is a value $m_r$ such that $f_r(x)=f_r^*(x)$ for all $x>m_r$.
 We further may assume that $r\mapsto m_r$ is a definable continuous and monotonic function. Thus we can find an interval $[a, b]$ and $m\in k$ 
 such that if $r\in[a,b]$, then $f_r(x)=f_r^*(x)$ for $x>\max(m,\alpha)$. Replacing $h_0$ by $f_a$, $h_1$ by $f_b$, and
 $\alpha$ by $\max(\alpha, m)$,  we may assume that $F$ fixes the graph $f_r$ for all $r$. 
 
 \medskip
Let $\psi:(\alpha,+\infty)\rightarrow (a,b)$ be an increasing $k$-definable continuous bijection.    
Define $$f(x)=f_{\psi(x)}(x) \hbox{ for } x>\alpha .$$   
Observe that for each $r$ the graph of $f_r$ intersects the graph of $f$ at a unique point $x_r$. 
Hence, for sufficiently large $x$, we have $f_{\psi(x)} (x)=f(x)$. Since for all $r$ we have $f_r(x)=f_r^*(x)$ for $x>\max(m,\alpha)$,  
it follows that $F(x,f(x))=(x,f(x))$.  But we have assumed that the hypothesis for Case 2 does not hold, that is, 
$F(x,y)\ne (x,y)$ for all $(x,y)\in C$, and thus we reach a contradiction.
\hfill\qed.

\section{Remarks and Questions}
\begin{enumerate}\item Using quantifier elimination for real closed fields and the decidability of the real algebraic numbers, 
 our construction could be done effectively, producing a computable rigid non-Archimedean real closed field.
 
\item Let $T$ be the theory of an o-minimal expansion of a real closed field in a countable language.  Our construction works in this 
 setting to build a rigid non-Archimedean $\M\models T$ of o-minimal dimension two.
 
\item Presumably our method can be extended to produce rigid real closed fields of any finite transcendence degree greater than one.
Can one produce a rigid countable non-Archimedean real closed field of infinite transcendence degree?
 It  does not seem to  be
 useful in the infinite transcendence degree case as we can no longer easily diagonalize over all possible automorphisms.\footnote{While this paper was under review, Michael Lange \cite{lange} has shown that indeed the methods can be generalized to produce rigid real closed fields of transcendence degree $n$ for each $n$.  He also developed a novel diagonalization strategy to produce a rigid model of transcendence degree $\aleph_0$.}
 
 \item {One might be tempted to try producing   larger rigid real closed fields using the following idea.
 Let $K$ be a countable  real closed subfield of $\R$ of positive transcendence degree. The arguments above allow us to construct 
 $K\< x,y\> $ with no non-trivial automorphism fixing $K$.  But, while $K$ is itself rigid, automorphisms of $K\< x,y\> $ need not fix $K$.

To illustrate, let $B$ be a transcendence base for $K$,  $b\in B$ and $B_0=B\setminus\{b\}$.  Let $K_0$ be the real closure of $\Q(B_0)$.
 Note that $b$ is the unique realization of $\tp(b/K_0)$ in $K$.  Let $L_0=K_0\< x\> $ and $L= L_0\< y\> .$
 
  For a real closed field $F$, we say that $a\in \calR\setminus F$ realizes a {\em non-cut} over $F$
 if $F<|a|$ or there is $b\in F$ with $|a-b|<\epsilon$ for all $\epsilon\in F$ with $\epsilon>0$.  Otherwise, we say $a$ realizes a {\em cut} over $F$.
 
 Since $\tp(x/K_0)$ is a non-cut and $\tp(b/K_0)$ is a cut,  by \cite{marker} 3.3, $\tp(b/K_0)$ is not realized in $L_0$.
 By construction, $\tp(y/L_0)$ has a unique realization in $L$, but   $b+{1\over x}$ is a second realization of $\tp(b/L_0)$ in $L_0\< b\> $.
 Thus, by \cite{marker} 3.6, $\tp(b/L_0)$ is not realized in $L$.  It follows that $b$ and $b+{1\over x}$ realize the same type over $L$.
 Thus there is an automorphism of $K\< x,y\> =L\< b\> $ fixing $L$ and sending $b$ to $b+{1\over x}$.}
 
 \item In any non-Archimedean real closed field $K$ there is an natural valuation where the valuation ring is 
 $${\cal O}= \{x\in K: |x|<n\hbox{ for some }n\in \Nn\}.$$ The residue field is always isomorphic to a real closed subfield of $\R$ and the value group is divisible.
 In our field the residue field is the real algebraic numbers and 
 the value group is $\Q$.  Are there any restrictions on the residue field and value groups of a rigid non-Archimedean real closed field?
\end{enumerate}

\end{document}